\documentclass[reqno]{amsart}
\begin{document}
\title[zero-pressure gas dynamics system]
{Spherically Symmetric Solutions of Multi-dimensional Zero-Pressure Gas Dynamics 
System}

\author[Pal Choudhury, Joseph and Sahoo]
{Anupam Pal Choudhury, K.T.Joseph and Manas R. Sahoo}

\address{ Anupam Pal Choudhury, K.T.Joseph and Manas R. Sahoo \newline
TIFR Centre for Applicable Mathematics\\
Sharada Nagar, Chikkabommasandra, GKVK P.O.\\
Bangalore 560065, India}
\email{\newline anupam@math.tifrbng.res.in,
ktj@math.tifrbng.res.in,
sahoo@math.tifrbng.res.in}

\thanks{}

\begin{abstract}

In this article we find explicit formulae for spherically symmetric solutions of the multidimensional zero-pressure gas dynamics
system and its adhesion approximation. The asymptotic behaviour of the explicit solutions of the adhesion approximation has been studied.
The radial components of the velocity and density satisfy a simpler equation which enables us to get explicit formulae for different types of domains and study its
asymptotic behaviour. A class of solutions for the inviscid system with conditions
on the mass instead of conditions at origin has also been analyzed. 
 
\end{abstract}

\maketitle
\numberwithin{equation}{section}
\numberwithin{equation}{section}
\newtheorem{theorem}{Theorem}[section]
\newtheorem{remark}[theorem]{Remark}

\section{Introduction}	

A proper understanding of radially symmetric solutions can be considered as a fundamental step in studying the multi-dimensional
systems. Though apparently these solutions seem to be one-dimensional in nature, but the singularities at $r=0$ make the study in these cases
much more complicated and more intriguing than in the general one-dimensional cases (see \cite{c1,pls1,mc1} for illuminating discussions on this). 
This article is devoted to the study of explicit spherically symmetric solutions for the multi-dimensional zero-pressure gas dynamics system 
and the associated adhesion model.\\
The adhesion model 
\begin{equation}
\begin{aligned}
u_t + (u.\nabla)u &=\frac{\epsilon}{2} \Delta u,\\
\rho_t + \nabla.(\rho u)&=0.
\end{aligned}
\label{e1.1}
\end{equation}
is one of the analytical models proposed to describe the large-scale structure of the universe. Here $u,\rho$ are the velocity and the density of the particles respectively (see  \cite{g1} and 
references therein for the physical importance and analysis of solutions). 
This system and its inviscid counterpart, the  multi-dimensional zero-pressure gas dynamics system 
\begin{equation}
\begin{aligned}
u_t + (u.\nabla)u &=0,\\
\rho_t + \nabla.(\rho u)&=0
\end{aligned}
\label{e1.2}
\end{equation}
have been active field of research since they were introduced, see \cite{a1,g1,g2,s1,w1}.
However there is no clear analytical understanding of the solutions of these 
equations. The well-posedness theory and large time behaviour of the solutions still call for a better understanding.
The advantage of
\eqref{e1.2} over \eqref{e1.1} lies in the fact that it is simpler to describe the solution of \eqref{e1.2} in the region where 
classical smooth solutions exist. It is well-known that global smooth solutions for \eqref{e1.2} do not exist
in $\{ (x,t) : x  \in R^n, t>0\}$, even when the initial data
\begin{equation}
\begin{aligned}
u(x,0)&=u_0(x),\\
\rho(x,0)&=\rho_0(x)
\end{aligned}
\label{e1.3}
\end{equation}
are smooth with compact support. Therefore we need to interpret the global solutions in a weak sense.
But then weak solutions again are not unique and additional 
conditions need to be specified to select the physically interesting solution. One way of selecting the \textquoteleft physical\textquoteright  solution 
of \eqref{e1.2} is by taking the limit of solution of \eqref{e1.1} as $\epsilon$ goes to $0$. 

 In the context of 
the large scale structure formation, the fastest growing mode in the linear theory 
has decaying vorticity. So it is natural to seek potential solutions of the model equations \eqref{e1.1} and \eqref{e1.2}.
Then the velocity can be 
represented in terms of a velocity potential $\phi$, see \cite{w1}. In this case,
$u^\epsilon$, $\epsilon>0$ can be constructed  using the standard Hopf-Cole transformation, 
\cite{j3, w1}.
Then the continuity equation for $\rho^\epsilon$ is a linear equation 
with smooth bounded coefficients which can 
be solved using the method of characteristics. We use this method to
construct weak asymptotic solution of \eqref{e1.1} and \eqref{e1.3}.

For the inviscid system \eqref{e1.2}, spherical solutions of 
the form $u=(x/r) q$, $\rho(x)=\rho(r)$, $r=|x|$, was constructed for $n=3$
in \cite{jm1}. This is a special case of potential velocity. 
 It was shown that the radial components of velocity and density satisfy  
\begin{equation}
\begin{aligned}
q_t +q q_r &=0, \\
\rho_t +\frac{1}{r^2}(r^2 \rho q)_r &=0,\,\,\, r>0, 
\end{aligned}
\label{e1.4}
\end{equation}
where $r=|x|$. By the change of variable $p=r^2\rho$ the equation 
\eqref{e1.4} takes the form
\begin{equation}
\begin{aligned}
q_t +q q_r &=0, \\
p_t +(q p)_r &=0,\,\,\, r>0.
\end{aligned}
\label{e1.5}
\end{equation} 
Classical theory of hyperbolic conservation laws of Lax \cite{la1} is not applicable due to the
formation of $\delta$ waves. However the system \eqref{e1.5} is well studied in \cite{j2,j4,j6}.
Using explicit 
solution of this system, we obtain radial solution of \eqref{e1.2}, with
initial data at $t=0$ and with different behaviours at the origin $x=0$. It is clear that in the region of smoothness, $u$ and
$\rho$ constructed using \eqref{e1.5} is a solution of the multidimensional system \eqref{e1.1}.
As the solution is not smooth, it is not obvious that the constructed distribution is a weak solution of the 
multidimensional system \eqref{e1.2}. 

The structure of this paper is as described below.

In section 2, we derive the equations for \eqref{e1.1} and \eqref{e1.2} in the radial cases. For the corresponding equation of \eqref{e1.2} 
a linearization using Hopf-Cole transformation is also provided. In section 3, we recall the notion of weak asymptotic solutions from \cite{a1,s1} and
recollect the Rankine-Hugoniot conditions derived therein for the inviscid system \eqref{e1.2}. Section 4 is devoted to the construction of explicit weak asymptotic solutions
for \eqref{e1.2} using the adhesion approximation. In Section 5, we solve boundary value problems for radial adhesion model in different interesting domains  and study its asymptotic behaviour for 
space dimensions $n=2$ and $n=3$. In Section 6, the Rankine-Hugoniot conditions for the one-dimensional systems derived from \eqref{e1.2} for the radial cases
have been discussed. In the last section, we find explicit solution for the radial inviscid system, with conditions on the mass. 

\section{Corresponding equations in the radial case and their linearization} 
To study the radial solutions of \eqref{e1.1} with prescribed initial and boundary conditions,
we look for radial components of velocity  and density of the form 
\begin{equation}
u(x,t)=\frac{x}{r} q(r,t),\,\,\,\rho(x,t)=\rho(r,t),\,\,\,r=|x|.
\label{e2.1}
\end{equation}
We find the equations for $q$ and $\rho$ and show that the resulting system can be linearized
using the Hopf-Cole transformation. More precisely, we have the following

\begin{theorem}\label{thm2.1}
The equations \eqref{e1.1}, for radial components, transform into the following system:
\begin{equation}
\begin{aligned}
&q_t +q q_r =\frac{\epsilon}{2} [q_{rr}+\frac{(n-1)}{r}q_r-\frac{(n-1)}{r^2}q],\\
&\rho_t +r^{-(n-1)}(r^{(n-1)} \rho q)_r =0,\,\,\, r>0,\,\,\,t>0.
\end{aligned}
\label{e2.2} 
\end{equation}
The above system can further be linearized by the transformation
\begin{equation}
 q(r,t)= -\epsilon \frac{a_r}{a},\,\,\, \rho(r,t)=r^{-(n-1)}p(r,t)
\label{e2.3}
\end{equation}
to the system
\begin{equation}
 a_t=\frac{\epsilon}{2}[a_{rr}+\frac{(n-1)}{r} a_r],\,\,\, p_t+(pq)_r=0.
\label{e2.4}
\end{equation}
\end{theorem}
\begin{proof}
A simple computation shows that
\begin{equation}
(u_j)_{x_k}=\begin{cases} (\frac{1}{r})q +(\frac{x_j^2}{r^2}) q_r -(\frac{x_j^2}{r^3})q,\,\,\ \,\,\ k=j\\
(\frac{x_j x_k}{r^2}) q_r -(\frac{x_j x_k}{r^3}) q,\,\,\,\,\ k\neq j
\end{cases}  
\label{e2.5}
\end{equation}
and
\begin{equation}
(u_j)_{x_k x_k}=\begin{cases} (\frac{x_j^3}{r^3}) q_{rr} +(\frac{3x_j}{r^2}-\frac{3x_j^3}{r^4}) q_r 
+(\frac{3x_j^3}{r^5}-\frac{3x_j}{r^3}) q,\,\,\ \,\,\, k=j\\
 (\frac{x_j x_k^2}{r^3}) q_{rr} +(\frac{x_j}{r^2}-\frac{x_j x_k^2}{r^4}-\frac{2x_jx_k^2}{r^4}) q_r +
(\frac{3x_jx_k^2}{r^5}-\frac{x_j}{r^3})q,\,\,\, \,\,k\neq j
\end{cases}
\label{e2.6}
\end{equation}
From \eqref{e2.1},\eqref{e2.5} and \eqref{e2.6}, we have for $j=1,2,...n$
\[
\begin{aligned}
(u_j)_t + \sum_{k=1}^n u_k (u_j)_{x_k} &=(\frac{x_j}{r}) q_t +
(\frac{x_j}{r}) q [(\frac{1}{r}) q +(\frac{x_j^2}{r^2}) q_r -(\frac{x_j^2}{r^3}) q]\\
&+q\sum_{k\neq j}\frac{x_k}{r}[(\frac{x_jx_k}{r^2})q_r-(\frac{x_jx_k}{r^3})q]\\
&=\frac{x_j}{r}[q_t+q q_r] 
\end{aligned} 
\]
and
\[
\begin{aligned}
\frac{\epsilon}{2} \Delta u_j &=\frac{\epsilon}{2} \frac{x_j}{r}[(\frac{x_j^2}{r^2})q_{rr}+(\frac{3}{r}-\frac{3x_j^2}{r^3})q_r
+(\frac{3x_j^2}{r^4}-\frac{3}{r^2})q]\\
&+\sum_{k\neq j}[(\frac{x_k^2}{r^2})q_{rr}+(\frac{1}{r}-\frac{3x_k^2}{r^3})q_r +(\frac{3x_k^2}{r^4}-\frac{1}{r^2})q]\\
&=\frac{\epsilon}{2} \frac{x_j}{r}[q_{rr}+\frac{(n-1)}{r}q_r-\frac{(n-1)}{r^2}q].
\end{aligned} 
\]
Using these in \eqref{e1.1}, we get
\begin{equation} 
\frac{x_j}{r} \{q_t +q q_r\} = \frac{\epsilon}{2} \frac{x_j}{r}[q_{rr}+\frac{(n-1)}{r}q_r -\frac{(n-1)}{r^2}q],\,\,\,j=1,2,3,...n,
\label{e2.7}
\end{equation}
Also for the continuity equation, we have
\begin{equation}
\begin{aligned}
\rho_t +\sum_{k=1}^n(u_k \rho)_{x_k}&=\rho_t +\sum_{k=1}^n\{(\frac{1}{r} q 
+(\frac{x_k^2}{r^2}) q_r 
-(\frac{x_k^2}{r^3}) q)\rho + (\frac{x_k^2}{r^2}) q \rho_r \}\\
&=\rho_t +\frac{n}{r} q\rho +q_r \rho -\frac{1}{r} q \rho+ q \rho_r\\
&=\rho_t +(\frac{n-1}{r}) q\rho+(q\rho)_r\\
&=\rho_t + {r^{-(n-1)}}(r^{(n-1)} \rho q)_r.
\end{aligned}
\label{e2.8}
\end{equation}
From \eqref{e2.7} and \eqref{e2.8}, it follows 
that $(u,\rho)$ of the form \eqref{e2.1} is a solution of \eqref{e1.1} iff
$q$ and  $\rho$ satisfy the equations 
\begin{equation}
\begin{aligned}
&q_t +q q_r =\frac{\epsilon}{2} [q_{rr}+\frac{(n-1)}{r}q_r-\frac{(n-1)}{r^2}q],\\
&\rho_t +r^{-(n-1)}(r^{(n-1)} \rho q)_r =0.
\end{aligned} 
\label{e2.9}
\end{equation}
Let $p=r^{(n-1)} \rho$. Then the above system becomes 
\begin{equation}
\begin{aligned}
&q_t +q q_r =\frac{\epsilon}{2} [q_{rr}+\frac{(n-1)}{r}q_r-\frac{(n-1)}{r^2}q],\\
&p_t + (pq)_r =0, \,\,\, r>0,\,\,\,t>0.
\end{aligned}
\label{e2.10} 
\end{equation}
 Now the equation for the velocity component can be linearized using Hopf-Cole transformation. For this 
first we note that if $Q(r,t)$ is a solution of
\begin{equation}
\begin{aligned}
Q_t +\frac{1}{2}((Q_r)^2) =\frac{\epsilon}{2} [Q_{rr}+\frac{(n-1)}{r}Q_r]
\end{aligned}
\label{e2.11}
\end{equation}
then $q(r,t)=Q_r(r,t)$ is a solution of the first equation of \eqref{e2.10}.
Now let
\begin{equation*}
 Q(r,t)=-\epsilon \log a^\epsilon.
\end{equation*}
An easy calculation then shows that $Q$ satisfies \eqref{e2.11} iff $a^\epsilon$ satisfies
\begin{equation}
a_t =\frac{\epsilon}{2} [a_{rr}+\frac{(n-1)}{r}a_r].
\label{e2.12}
\end{equation}
\end{proof}

Proceeding as in the proof of Theorem $(2.1)$, we can prove the following result for the inviscid case.
\begin{theorem}
The equations \eqref{e1.2}, for radial components, transform into the system
\begin{equation}
\begin{aligned}
&q_t +q q_r =0,\\
&\rho_t +r^{-(n-1)}(r^{(n-1)} \rho q)_r =0.
\end{aligned} 
\label{e2.13}
\end{equation}
Using the transformation $\rho(r,t)=r^{-(n-1)}p(r,t)$, the second equation above can again be written as $$p_t+(pq)_r=0.$$
\end{theorem}

\section{Weak asymptotic solutions for the multi-dimensional case}
In this section, we recollect the notion of weak asymptotic solutions and generalized $\delta-$shock wave type solutions for the inviscid system
\eqref{e1.2} as introduced in \cite{a1,s1}.\\
\textbf{Definition $3.1$} (Weak asymptotic solutions)\\
A family of smooth functions $(u^{\epsilon},\rho^{\epsilon})_{\epsilon >0}$ is called a \textit{weak asymptotic solution} of the system
\eqref{e1.2} with initial conditions \eqref{e1.3} if  
\begin{equation}
\begin{aligned}
u^\epsilon_t + (u^\epsilon.\nabla)u^\epsilon &=o_{D'}(R^n)(1),\\
\rho^\epsilon_t + \nabla.(\rho^\epsilon u^\epsilon) &=o_{D'}(R^n)(1),\\
u^\epsilon(x,0)-u_0(x) &=o_{D'}(R^n)(1),\\
\rho^\epsilon(x,0)-\rho_0(x) &=o_{D'}(R^n)(1)
\end{aligned}
\label{e3.1}
\end{equation}
where $\langle o_D'(R^n)(1),\eta\rangle \rightarrow 0$ as $\epsilon \rightarrow 0$
uniformly in $t>0$, for every $\eta \in C_{c}^\infty(R^n)$.\\
\\
Now let us suppose that the initial conditions \eqref{e1.3} are of the form
\begin{equation}
\begin{aligned}
 u_{0}(x)&=u_{+}(x)H(S(x,0)>0)+u_{-}(x)H(S(x,0)<0),\\
\rho_{0}(x)&=\bar{\rho}_{+}(x)H(S(x,0)>0)+\bar{\rho}_{-}(x)H(S(x,0)<0)+ \hat{e}\delta_{S(x,0)=0}
\end{aligned}
\label{e3.2}
\end{equation}
where $u_{+}$,$u_{-}$,$\bar{\rho}_{+}$,$\bar{\rho}_{-}$ are smooth functions away from the surface $S(x,0)=0$. 
We consider an ansatz for \eqref{e1.2} of the form
\begin{equation}
\begin{aligned}
  u(x,t)&=u_{+}(x,t)H(S(x,t)>0)+u_{-}(x,t)H(S(x,t)<0),\\
 \rho(x,t)&=\bar{\rho}_{+}(x,t)H(S(x,t)>0)+\bar{\rho}_{-}(x,t)H(S(x,t)<0)+\hat{e}(t)\delta_{S(x,t)=0}
\end{aligned}
\label{e3.3}
\end{equation}
where $u_{+}$,$u_{-}$,$\bar{\rho}_{+}$ and $\bar{\rho}_{-}$ are smooth functions away from the surface $S(x,t)=0$.\\
\\
\textbf{Definition $3.2$} (Generalized $\delta-$shock wave type solution)\\
The ansatz \eqref{e3.3} is called a \textit{generalized $\delta-$shock wave type solution} for the Cauchy problem \eqref{e1.2},\eqref{e1.3}
if it is the distributional limit of a weak asymptotic solution $(u^{\epsilon},\rho^{\epsilon})_{\epsilon >0}$ as $\epsilon \rightarrow 0$.\\
\\
In \cite{a1} it was further remarked that a \textit{generalized $\delta-$shock wave type solution} for \eqref{e1.2},\eqref{e1.3} can also 
be defined in terms of the following integral identities:\\
Let $\Gamma_t =\{(x,t) : S(x,t)=0\}$ be a surface of discontinuity, $\sigma(x,t)=\frac{dS(x,t)}{|\nabla_{x,t}S|}$ be the
surface measure and $\frac{D}{Dt}=\partial_t + \frac{u_{+}+u_{-}}{2}.{\nabla_{x}}$, $u_{-}$ and $u_{+}$ are velocity behind and
ahead of the discontinuity. Then the ansatz \eqref{e3.3} is a generalized $\delta-$shock wave type solution of \eqref{e1.2},\eqref{e1.3} if the 
following integral identities hold:
\begin{equation}
 \begin{aligned}
&\int_{\Omega\times[0,T]-{\Gamma_t}}(\rho \phi_t+\rho u.\nabla_x(\phi))dx dt+\int_{\Gamma_t}\hat{e}(t) \frac{D\phi}{Dt}\sigma(x,t)=0,\\
&\int_{\Omega\times[0,T]-{\Gamma_t}}(u_t+(u.\nabla)u)\phi dx dt-\int_{\Gamma_t} [u]\frac{DS}{Dt}\phi \sigma(x,t)=0
\end{aligned}
\label{e3.4}
\end{equation}
for all $\phi \in C_c^{\infty}(\Omega \times (0,T))$.\\
\\
In \cite{a1}, using the work of Majda \cite{ma1,ma2} on existence and stability of multidimensional shock fronts, the authors gave a local construction of $(u,\rho)$ and a surface $S(x,t)=0$, from an initial discontinuous wave front \eqref{e3.2} satisfying an entropy condition. The solution satisfied \eqref{e1.2} in the region of smoothness and 
the Rankine Hugoniot conditions
\begin{equation}
\begin{aligned}
 S_t + u_{\delta}.\nabla_{x}S|_{\Gamma_t}&=0,\\
\frac{\delta \hat{e}}{\delta t} + \nabla_{\Gamma_t}(\hat{e} u_{\delta})&=([\bar{\rho} u]-[\bar{\rho}]u_{\delta}).\nabla_{x}S|_{\Gamma_t},
\end{aligned}
\label{e3.5}
\end{equation}
along the surface of discontinuity $S(x,t)=0$.\\
In Section $7$, we construct radial solutions which satisfy  
\eqref{e1.2} in the region of smoothness and the Rankine-Hugoniot conditions \eqref{e3.5}
along the surface of discontinuity.

\section{Explicit weak asymptotic solutions for the multidimensional case using adhesion approximation}
In this section, we construct explicit weak asymptotic solution of 
the system \eqref{e1.2} with initial conditions \eqref{e1.3} using
the adhesion approximation \eqref{e1.1}, under the additional condition that
$u=\nabla_{x}\phi$.

We recall from Section 3 that a weak asymptotic solution for the system \eqref{e1.2} with initial conditions \eqref{e1.3} is a family of smooth functions 
$(u^\epsilon,\rho^\epsilon)_{\epsilon>0}$
satisfying the conditions \eqref{e3.1}. We consider the adhesion approximation \eqref{e1.1} 
with initial conditions which are regularizations of \eqref{e1.3}.
The velocity component $u^\epsilon$ can be constructed   
using Hopf-Cole transformation, and then the continuity equation for $\rho^\epsilon$ is 
solved using the method of characteristics. The resulting
family $(u^\epsilon,\rho^\epsilon)_{\epsilon>0}$, is a weak asymptotic solution of \eqref{e1.2} and 
\eqref{e1.3} as proved in the next theorem. 

\begin{theorem}\label{thm4.1}

Assume $u_0(x)=\nabla_{x}\phi_0$ where $\phi_0 \in W^{1,\infty}(R^n)$  and $\rho_0 \in L^{\infty}(R^n)$. Let
$\phi_0^\epsilon =\phi_0*\eta^\epsilon$, $\nabla_{x}\phi_0^\epsilon 
=\nabla_{x}(\phi_0*\eta^\epsilon)$ and 
$\rho_0=\rho_0*\eta^\epsilon$, where $\eta^\epsilon$ is the usual Friedrichs 
mollifier in the space variable $x \in R^n$.\ 
Further let

\begin{equation}
\begin{aligned}
u^\epsilon(x,t)&=\frac{\int_{R^n}(\nabla_{y}\phi_0^\epsilon(y))e^{\frac{-1}{\epsilon}\left[\frac{|x-y|^2}{2t}+\phi_0^\epsilon(y)\right]}dy}
{\int_{R^n}e^{\frac{-1}{\epsilon}\left[\frac{|x-y|^2}{2t}+\phi_0^\epsilon(y)\right]}dy},\\
\rho^\epsilon(x,t)&=\rho_0^\epsilon(X^\epsilon(x,t,0)) J^\epsilon(x,t,0),
\end{aligned}
\label{e4.1}
\end{equation}
where $X^\epsilon(x,t,s)$ is the solution of 
$\frac{d X^{\epsilon}(s)}{ds}= u^\epsilon(X^{\epsilon},s)$ with $X^{\epsilon}(s=t)=x$ and  
$J^\epsilon(x,t,0)$ is the Jacobian matrix of $X^\epsilon(x,t,0)$ w.r.t. 
$x$. 

Then 
$(u^\epsilon,\rho^\epsilon)$ is a weak asymptotic 
solution to \eqref{e1.2} and \eqref{e1.3}.
\end{theorem}

\begin{proof}
Following \cite{h1,j3,w1}, we note that if $\phi^\epsilon$ is a solution of 

\begin{equation}
\phi_t +\frac{|\nabla_{x} \phi|^2}{2}=\frac{\epsilon}{2} \Delta \phi,\,\,
\phi(x,0)=\phi_0^\epsilon(x),
\label{e4.2}
\end{equation}
then $u^\epsilon =\nabla_{x}\phi^\epsilon$ is a solution of \eqref{e1.1} with initial condition
$u^\epsilon(x,0) =\nabla_{x}\phi_0^\epsilon(x)$.
Now using the Hopf-Cole transformation $\theta =e^{-\frac{\phi}{\epsilon}}$, 
it follows that $\phi^\epsilon$ is the solution of \eqref{e4.2} iff $\theta^\epsilon$ is the 
solution of
\begin{equation}
\theta_t =\frac{\epsilon}{2}\Delta \theta, \,\,\,
\theta(x,0)=e^{-\frac{\phi_0^\epsilon(x)}{\epsilon}}.
\label{e4.3}
\end{equation}
Solving \eqref{e4.3}, we get,
\begin{equation}
\theta^\epsilon(x,t)=\frac{1}{(2\pi t\epsilon)^{n/2}}\int_{R^n} 
e^{-\frac{1}{\epsilon}[\frac{|x-y|^2}{2t}+ \phi_0^\epsilon(y)]} dy
\label{e4.4}
\end{equation}

Now
\begin{equation}
\begin{aligned}
\theta^{\epsilon}(x,t)_{x_j}&=\frac{1}{(2\pi t\epsilon)^{n/2}}\int_{R^n} 
\partial_{x_j}(e^{-\frac{1}{\epsilon}[\frac{|x-y|^2}{2t} 
+\phi_0^\epsilon(y)]}) dy\\
&=-\frac{1}{(2\pi t\epsilon)^{n/2}}\int_{R^n} 
(\partial_{y_j}(e^{-\frac{1}{\epsilon}\frac{|x-y|^2}{2t}})).e^{-\frac{1}{\epsilon} 
\phi_0^\epsilon(y)} dy\\
&=\frac{1}{(2\pi t\epsilon)^{n/2}}\int_{R^n} 
-(\frac{1}{\epsilon})(\partial_{y_j} \phi_0^\epsilon(y))
e^{-\frac{1}{\epsilon}\frac{|x-y|^2}{2t}}.e^{-\frac{1}{\epsilon} 
\phi_0^\epsilon(y)} dy.
\label{e4.5}
\end{aligned}
\end{equation}

In the last line we have used integration by parts with respect to the variable $y$. Since $\theta =e^{-\frac{\phi}{\epsilon}}$, we get from 
\eqref{e4.4} and \eqref{e4.5},

\begin{equation}
u^\epsilon(x,t)=\frac{\int_{R^n}(\nabla_{y}\phi_0^\epsilon(y))e^{\frac{-1}{\epsilon}\left[\frac{|x-y|^2}{2t}+\phi_0^\epsilon(y)\right]}dy}
{\int_{R^n}e^{\frac{-1}{\epsilon}\left[\frac{|x-y|^2}{2t}+\phi_0^\epsilon(y)\right]}dy}.
\label{e4.6}
\end{equation}

From this formula it is clear that

\begin{equation}
||u^\epsilon||_{L^\infty(R^n\times[0,\infty))}\leq ||\nabla_{x} 
\phi_0^\epsilon||_{L^\infty(R^n)}.
\label{e4.7}
\end{equation}

Further as differentiation under the integral sign is 
justified if $\nabla_x\phi_0$ is bounded measurable, it follows that $u^\epsilon$ is a $C^\infty$ function in 
$R^n\times(0,\infty)$. Indeed, we can write the formula 
\eqref{e4.6} for $u^\epsilon$ in the following form

\begin{equation}
u^\epsilon(x,t)=\frac{\int_{R^n}\nabla_{x}\phi_0^\epsilon(x-\sqrt{(2 t)}y)e^{\frac{-1}{\epsilon}\left[|y|^2+\phi_0^\epsilon(x-\sqrt{(2 
t)}y)\right]}dy}
{\int_{R^n}e^{\frac{-1}{\epsilon}\left[|y|^2+\phi_0^\epsilon(x-\sqrt{(2 t) }y)\right]}dy}.
\label{e4.8}
\end{equation}

Differentiating \eqref{e4.8}, using chain rule, it follows that

\begin{equation}
|\partial_x^\alpha \partial_t^j u^\epsilon| \leq 
\frac{C_{\alpha,j}}{\epsilon^{|\alpha|+j}},
\label{e4.9}
\end{equation}
where $C_{\alpha,j}$ depends only on $||\nabla_{x} \phi_0^\epsilon||_{L^\infty}$.
Next we consider the continuity equation with coefficient $u^\epsilon$ and 
regularized initial condition,

\begin{equation}
\rho_t^\epsilon +\nabla_x. (u^\epsilon 
\rho^\epsilon)=0,\,\,\,\rho^\epsilon(x,0)=\rho_0^\epsilon(x).
\label{e4.10}
\end{equation}

As $u^\epsilon$ is smooth, we use the method of characteristics to find $\rho$. Let 
$X^{\epsilon}(x,t,s)$ be the solution of 

\begin{equation}
\frac{d X^{\epsilon}(s)}{ds}= u^\epsilon(X^{\epsilon},s)\,\,\, , X^{\epsilon}(s=t)=x.
\label{e4.11}
\end{equation}

Since $u^\epsilon$ satisfies the estimates \eqref{e4.7} and \eqref{e4.9},
by the existence and uniqueness theory of ODE, there exists a unique 
solution $X^{\epsilon}(s)$ to \eqref{e4.11} for all $0\leq s \leq t$. To make 
the dependence of $x$ and $t$ 
explicit, let us denote it by $X^{\epsilon}(x,t,s)$. This flow takes the point $(x,t)$ to
an initial point $(X^{\epsilon}(x,t,0),0)$ and conversely. 

Let $J(X^\epsilon(x,t,s))$ be the Jacobian determinant of 
$X^\epsilon(x,t,s)$ with respect to $x$. Then

\begin{equation}
\rho^\epsilon(x,t)= \rho_0^\epsilon(X^{\epsilon}(x,t,0)) J(X^\epsilon(x,t,0))
\label{e4.12}
\end{equation}
is the solution of \eqref{e4.10}.\\
The family $(u^\epsilon,\rho^\epsilon)_{\epsilon>0}$ given by \eqref{e4.6} and
\eqref{e4.12} is then a weak 
asymptotic solution.This follows easily as 
\[
\epsilon \int_{R^n}\Delta u^\epsilon \eta(x) dx =
\epsilon \int_{R^n} u^\epsilon \Delta \eta(x) dx =O(1)\epsilon 
\]
uniformly in $t$ for every $\eta \in C_0^\infty(R^n)$ as $u^\epsilon$
is bounded independent of $\epsilon$ and $(x,t)$ by the estimate \eqref{e4.7}.
The solution $(u^\epsilon,\rho^\epsilon)_{\epsilon>0}$ satisfies
the initial conditions since we have
\begin{equation*}
\begin{aligned}
u^\epsilon(x,0)-\nabla_{x}\phi_0(x)&=u^\epsilon(x,0)-\nabla_{x}\phi_0^\epsilon(x) 
+\nabla_{x}\phi_0^\epsilon(x) -\nabla_{x}\phi_0(x)\\
&=\nabla_{x}\phi_0^\epsilon(x) -\nabla_{x}\phi_0(x),\\
\rho^\epsilon(x,0)-\rho_0(x)&=\rho^\epsilon(x,0)-\rho_0^\epsilon(x)+
\rho_0^\epsilon(x)-\rho_0(x)=\rho_0^\epsilon(x)-\rho_0(x)
\end{aligned}
\end{equation*}
and $\nabla_{x}\phi_0^\epsilon(x) -\nabla_{x}\phi_0(x)$, 
$\rho_0^\epsilon(x)-\rho_0(x)$ go to zero in distributions.
\end{proof}
There is a large class of interesting initial data  which admit gradient type solutions
as described in the theorem. One such case is when the initial data is of radial type.

\begin{equation}
 u(x,0)=u_0(x)=\frac{x}{r}q_0(r),
\label{e4.13}
\end{equation}
Clearly it can be written as a gradient,
\[
 u_0(x)=\nabla_{x}\phi_0(x)=\nabla_{x}(\int_0^{|x|} q_0(s) ds).
\]
For density we give the initial condition
\begin{equation}
 \rho(x,0)=\rho_0(x).
\label{e4.14}
\end{equation}

\begin{theorem}\label{theorem4.2}
 The solution of \eqref{e1.1}, with initial conditions \eqref{e4.13} and \eqref{e4.14} with
$\int_0^{\infty} q_0(s) ds<\infty$ and $\int_{R^n}\rho_0(x) dx <\infty$ has the 
following asymptotic behaviour. The velocity component goes to $0$ as $t$ tends to $\infty$
uniformly on compact subsets of $R^n$ and the mass $\int_{R^n} \rho(x,t) dx$ is conserved.
\end{theorem}
\begin{proof}
From \eqref{e4.4}-\eqref{e4.5} and $\theta =e^{-\frac{\phi}{\epsilon}}$, we get the following formula for the velocity 

\begin{equation*}
u^\epsilon(x,t)=\frac{\int_{R^n}(\frac{x-y}{t}) e^{\frac{-1}{\epsilon}\left[\frac{|x-y|^2}{2t}+ \int_0^{|y|} q_0(s) ds \right]}dy}
{\int_{R^n}e^{\frac{-1}{\epsilon}\left[\frac{|x-y|^2}{2t}+\int_0^{|y|} q_0(s) ds \right]}dy}.
\end{equation*}
Making the change of variable $y\rightarrow (x-y)/\sqrt{2t}$, we get,

\begin{equation}
u^\epsilon(x,t)=\frac{\sqrt{\frac{2}{t}}\int_{R^n} y e^{\frac{-1}{\epsilon}\left[|y|^2+\int_0^{|x-\sqrt{(2 
t)}y|} q_0(s) ds\right]}dy}
{\int_{R^n}e^{\frac{-1}{\epsilon}\left[|y|^2+\int_0^{|x-\sqrt{(2 
t) }y|} q_0(s) ds\right]}dy}
\label{e4.15}
\end{equation}
From \eqref{e4.15}, it follows that,
as $t$ goes to $\infty$, $u^\epsilon$ goes to $0$ uniformly on compact subsets of $R^n$ . 
From the formula \eqref{e4.12} we have
\begin{equation*}
 \begin{aligned}
  \int_{R^n}\rho^{\epsilon}(x,t)dx&=\int_{R^n}\rho_0^\epsilon(X^{\epsilon}(x,t,0)) J(X^\epsilon(x,t,0))\\
&=\int_{R^n}\rho_0(x)dx, 
 \end{aligned}
\end{equation*}
 since for every $t$, $X^{\epsilon}(x,t,0):R^n \rightarrow R^n$ 
is a diffeomorphism and the determinant of Jacobian, $J(X^\epsilon(x,t,0))$ is positive.
 Thus the mass $\int_{R^n} \rho(x,t)dx$ is conserved.
\end{proof}

From the above discussions, we know that solution of the adhesion approximation \eqref{e1.1} has an explicit formula
\begin{equation}
\begin{aligned}
u^\epsilon(x,t)&=\frac{\int_{R^n}(\frac{x-y}{t})e^{\frac{-1}{\epsilon}\left[\frac{|x-y|^2}{2t}+\phi_0^\epsilon(y)\right]}dy}
{\int_{R^n}e^{\frac{-1}{\epsilon}\left[\frac{|x-y|^2}{2t}+\phi_0^\epsilon(y)\right]}dy},\\
\rho^\epsilon(x,t)&=\rho_0^\epsilon(X^\epsilon(x,t,0)) J^\epsilon(x,t,0),
\end{aligned}
\label{e4.16}
\end{equation}
We conclude this section with a result concerning the vanishing viscosity limit of the velocity component $u$.\\
For the general data $u(x,0)=\nabla_x \phi_0(x)$, formula \eqref{e4.16} implies 
\begin{equation}
 \lim_{\epsilon \rightarrow 0} u^\epsilon(x,t) =u(x,t)=\frac{x-y(x,t)}{t}
\notag
\end{equation}
where $y(x,t)$ is a minimizer in
\begin{equation}
\min_{y \in R^n} \{\phi_0(y)+\frac{|x-y|^2}{2t}\}.
\notag
\end{equation}
For almost every $(x,t)$ this minimizer is unique and so $u(x,t)$ is well defined a.e.
We show that if the initial data is radial then the solution is also radial and has a simpler form. 
\begin{theorem}\label{thm4.3} Assume the initial data $\phi_0(x)=\phi_0(|x|)$, then 
\begin{equation}
 \min_{y \in R^n} \{\phi_0(|y|)+\frac{|x-y|^2}{2t}\}=\min_{r\geq 0}\{\phi_0(r)+\frac{(|x|-r)^2}{2t}\}
\notag 
\end{equation}

\begin{equation}
u(x,t)=\frac{x-y(x,t)}{t} = (\frac{|x|-r(|x|,t)}{t}) \frac{x}{|x|}
\end{equation}
where $y(x,t)$ and $r(|x|,t)$ are the minimizers on the left and right of the first equality above, respectively.
\end{theorem}

\begin{proof}
 Since $||x|-|y||\leq |x-y|$, we have  
\begin{equation}
\frac{||x|-|y||^2}{2t}+\phi_0(y)\leq \frac{{|x-y|}^2}{2t} + \phi_0(|y|)
\notag
\end{equation}
and so
\begin{equation}
\min _{y\in R^n}\frac{||x|-|y||^2}{2t}+\phi_0(y)\leq \min_{r\in R^n}\frac{{|x-y|}^2}{2t} + \phi_0(|y|)
\notag
\end{equation}
Now consider $y=r \frac{x}{|x|}$, then $|x-y|=|\frac{x}{|x|}(|x|-r)|=||x|-|y||$ and so
\begin{equation}
\min_{y \in R^n} \{\phi_0(|y|)+\frac{|x-y|^2}{2t}\}\leq\min_{r\geq 0}\{\phi_0(r)+\frac{(|x|-r)^2}{2t}\} 
\notag
\end{equation}
and the first equality follows. Further the minimum is achieved on $y(x,t)$ of the form $y(x,t)=\frac{x}{|x|}r(|x|,t)$. Therefore
we have
\begin{equation}
(x-y(x,t))=x-\frac{x}{|x|}r(|x|,t) = \frac{x}{|x|}(|x|-r(|x|,t))
\notag
\end{equation}
and the second equality follows.
\end{proof}

{\bf Remark :} The weak limit of $\rho^\epsilon$, as $\epsilon$ goes to zero is a measure 
and passage to the limit in $u^\epsilon.\rho^\epsilon$ is still open. For the radial case considered in Section 7, we construct
a solution using the velocity component $u$ obtained above in theorem 4.3. 

\section{Radial solutions of the adhesion model with boundary conditions} 
In this section, we obtain explicit solutions for \eqref{e1.1} in some spherically
symmetric domains of the form $\Omega = D\times (0,\infty)$ for $2$ and $3$ space dimensions.\\

\textbf{Case 1:} Let $\Omega = \{(x,t) : |x| = r <R, t>0 \}$ be the domain under consideration.
We take initial conditions of the form
\begin{equation}
\begin{aligned}
  u(x,0)&=\frac{x}{|x|} q_0(r),\\
 \rho(x,0)&=\rho_0(r) ,\ \ \ \ |x|< R,
\end{aligned}
\label{e5.1}
\end{equation}
and boundary conditions 
\begin{equation}
\begin{aligned}
 \lim_{|x|\rightarrow R} \frac{x}{|x|}u(x,t)&=q_B,\ \text{where} \ q_B \ \text{is constant}\\
 \lim_{|x|\rightarrow R}\rho(x,t)&=\rho_B(t),\,\ if \ q_B<0.
\end{aligned}
\label{e5.2}
\end{equation}
(Here we note that a boundary condition for $\rho$ needs to be specified only if $q_B$ is negative.)\\
We also assume a first order consistency condition of the initial and boundary data at $\{(x,t): |x|=R, t=0 \}$.\\
Let $\Gamma = \{(x,t) \in \Omega : \textnormal{either}\,\,\, |x|=R \,\,\ \textnormal{or}\,\,\, t=0\}$. Define $p_{\Gamma}$ on $\Gamma$ by
\begin{equation}
p_{\Gamma}(r,t) = \begin{cases}
 p_0(r),&\text{if } t=0,\\
  p_B(t),&\text{if } r=R. 
\end{cases}
\label{e5.3}
\end{equation}
As $q$ is smooth, the characteristic $\beta(r,t,s)$ of $p_t +(qp)_r =0$ passing through $(r,t)$ exists for all $s<t$, till it meets $\Gamma$, at time $t=t_0(r,t) \geq 0$.\\
Let $J_k$, for  $k=0,1$, denote the Bessel functions of first kind and order k. With these notations, we have the following theorem. 
\begin{theorem}
Explicit solution for \eqref{e1.1}, with initial and boundary conditions \eqref{e5.1} - \eqref{e5.2} is given by
\begin{equation}
\begin{aligned}
 u(x,t)&=-\epsilon \frac{x}{|x|}\frac{\int_0^R \partial_r{G(r,\xi,t)}e^{-\frac{1}{\epsilon}\int_0 ^{\xi} q_0(s)ds}d\xi}
{\int_0^R G(r,\xi,t) e^{-\frac{1}{\epsilon}\int_0 ^\xi q_0 (s)ds} d\xi}\\
\rho(x,t)&=r^{-(n-1)} p_{\Gamma}(\beta(r,t,t_0 (r,t)),t_0(r,t))e^{-{\int_{t_{0}(r,t)}^t q_r(\beta(r,t,s)) ds}},
\label{e5.4}
\end{aligned}
\end{equation}
where $G(r,\xi,t)$ is as described below.\\
For $n=2$, $G$ has the form
\begin{equation}
G(r,\xi, t)=\frac{2}{R^2}\displaystyle {\sum_{n=1} ^{\infty}} \frac{\mu_n ^2 \xi }{((q_B/\epsilon)^2 R^2 + \mu_n ^2)J_0 ^2(\mu _n)}
J_0 (\frac {\mu_n r} {R})J_0 (\frac {\mu_n \xi} {R})e^{-\frac{\epsilon \mu_n ^2 t}{2 R^2}},
\label{e5.5}
\end{equation}
where $\mu_n$  are positive solutions of the transcendental  equation 
$\mu J_1 (\mu)-kR J_0(\mu)=0$ and for $n=3$, $G$ takes the form

\begin{equation}
G(r,\xi, t)=\frac{2\xi}{R r}\displaystyle {\sum_{n=1} ^{\infty}} \frac{\mu_n ^2 +((q_B/\epsilon)R-1)^2 }{\mu_n ^2 + (q_B/\epsilon)R((q_B/\epsilon)R -1)}
\sin (\frac {\mu_n r} {R}) \sin (\frac {\mu_n \xi} {R})e^{-\frac{\epsilon \mu_n ^2 t}{2R^2}},
\label{e5.6}
\end{equation}
where $\mu_n$  are positive solutions of the transcendental  equation 
$\mu \cot(\mu)+(q_B/\epsilon)R-1=0$.
\end{theorem}
\begin{proof}
By Theorem 2.1, it follows that the velocity $u$ is given by
$u(x,t)=-\epsilon \frac{x}{|x|}\frac{a_r}{a}$,
where $a$ satisfies the linear problem
\begin{equation}
\begin{aligned}
&a_t = \frac{\epsilon}{2}[a_{rr} +\frac{(n-1)}{r}a_r],\,\,\, r<R,\,\,\,t>0\\
&a(r,0)=e^{-\frac{\int_0^ r q(s) ds}{\epsilon}}, \,\,\,r<R\\
&\epsilon a_r(R ,t)+q_B a(R ,t)=0, \,\,\,t>0.\\
\end{aligned}
\label{e5.7} 
\end{equation}
Then $G$ given by \eqref{e5.5} and \eqref{e5.6} are just the Green's functions for the boundary value
problem \eqref{e5.7}, for $n=2$ and $n=3$ respectively, see \cite{p1}.  The formula \eqref{e5.4} for $u$ follows. 

Now we find a formula for $\rho$. Notice that $\rho(x,t)=\frac{1}{r^{(n-1)}}p$ where $p$ satisfies the equation
\[
p_t+qp_r=-q_r p.
\]
Using the method of characteristics we have $\frac{dp}{ds} =-q_r p$ along the curve $\frac{d\beta}{ds} =q$. Let us consider the characteristic curve $(\beta(r,t,s),s)$ passing through $(x,t)$. 
Integrating the equation for $p$, along this characteristic, we get $$p(r,t)=p(\beta(r,t,t_0)) e^{-\int_{t_0}^t q_{r}(\beta(r,t,s))ds},$$ where
$t_0=t_0(r,t)(<t)$ is the time  the curve touches the boundary $\Gamma$.
 \end{proof}

\textbf{Case 2:} Next let us consider the equation \eqref{e1.1} in the domain $\Omega=\{(x,t): R_1 < |x|< R_2,\,\,\, t>0\}$
with initial conditions 
\begin{equation}
\begin{aligned}
  u(x,0)&=\frac{x}{|x|} q_0(r),\\ 
 \rho(x,0)&=\rho_0(r) , \ \ R_1<|x|< R_2
\end{aligned}
\label{e5.8}
\end{equation}
and boundary conditions 
\begin{equation}
 \lim_{|x|\rightarrow R_i} \frac{x}{|x|}u(x,t)=q_i,\ i=1,2
\label{e5.9}
\end{equation}
where $q_i$ are constants. For $\rho$ we need to prescribe a boundary condition on $(|x|=R_i)$ only if $(-1)^{i+1}q_i$ is positive,
\begin{equation}
 \lim_{|x|\rightarrow R_i}\rho(x,t)=\rho_i (t),\,\,\,if\ (-1)^{i+1}q_i>0.
\label{e5.10}
\end{equation}
We also assume a first order consistency condition of the initial and boundary data at $\{(x,t): |x|=R_1,R_2 \, \textnormal{and} \, t=0 \}$.\\
Let $\Gamma = \{(x,t) \in \Omega : \textnormal{either}\,\,\, |x|=R_i \,\,\ \textnormal{or}\,\,\, t=0\}$. Define $p_{\Gamma}$ on $\Gamma$ by
\begin{equation}
p_{\Gamma}(r,t) = \begin{cases}
 p_0(r),&\text{if } t=0,\\
  p_i(t),&\text{if } r=R_i. 
\end{cases}
\label{e5.11}
\end{equation}
As $q$ is smooth,
let $\beta(r,t,s)$ be the characteristic of $p_t +(qp)_r =0$ passing through $(r,t)$ and let this characteristic meet the boundary
point at time $t=t_0(r,t)$.\\
Let $J_k$ ($k=0,1$) denote the Bessel functions of first kind and order $k$ and $Y_k$ ($k=0,1$) denote the Bessel functions of second
kind and order $k$. With these notations, we have the following theorem.
\begin{theorem}
Explicit solution for \eqref{e1.1}, with initial and boundary conditions \eqref{e5.1} - \eqref{e5.2} is given by
\begin{equation}
\begin{aligned}
 &u(x,t)=-\epsilon \frac{x}{|x|}\frac{\int_{R_1}^{R_2} \partial_r{G(r,\xi,t)} e^{-\frac{1}{\epsilon}\int_0 ^\xi q_0 (s)ds}d\xi}
{\int_{R_1}^{R_2} G(r,\xi,t) e^{-\frac{1}{\epsilon}\int_0 ^\xi q_0 (s)ds} d\xi}\\
&\rho(x,t)=r^{-(n-1)}p_{\Gamma}(\beta(r,t,t_0 (r,t)),t_0(r,t)) e^{-{\int_{t_0(r,t)}^t q_r(\beta(r,t,s)) ds}},
\label{e5.12}
\end{aligned}
\end{equation}
where $G(r,\xi,t)$ is as described below.\\
For $n=2$, $G$ has the form
\begin{equation}
 G(r,\xi,t)=\frac{\pi^2}{2}\sum_{n=1}^{\infty}\frac{\lambda_n^2}{B_n}[(q_2/\epsilon) J_0(\lambda_n r_2) - \lambda_n J_1(\lambda_n r_2)]^2 
\xi H_n(r) H_n(\xi) e^{-\frac{\epsilon \lambda_n^2 t}{2}}
\label{e5.13}
\end{equation}
where
\begin{equation}
\begin{aligned}
 B_n &= (\lambda_n ^2  +(q_2/\epsilon)^2)[-(k_1/\epsilon) J_0 (\lambda_n R_1)+\lambda_n J_1(\lambda_n R_1)]^2 \\
&-(\lambda_n^2+(k_1/\epsilon) ^2)
[(q_2/\epsilon) J_0 (\lambda_n R_2)-\lambda_n J_1(\lambda_n R_2)]^2,\\
H_n(r) &=[-(q_1/\epsilon) Y_0 (\lambda_n R_1)+\lambda_n Y_1 (\lambda_n R_1)]J_0 (\lambda_n r)\\
&-[-(q_1/\epsilon) J_0 (\lambda_n R_1)+
\lambda_n J_1 (\lambda_n R_1)]Y_0 (\lambda_n r)
\end{aligned}
\label{e5.14}
\end{equation}
and the $\lambda_n$ are positive roots of the equation
\begin{equation}
\begin{aligned}
 &[-(q_1/\epsilon) J_0(\lambda R_1)+\lambda J_1(\lambda R_1)][(q_2/\epsilon) Y_0(\lambda R_2)-\lambda Y_1(\lambda R_2)]\\
&-[(q_2/\epsilon) J_0(\lambda R_2)-\lambda J_1(\lambda R_2)][-(q_1/\epsilon) Y_0(\lambda R_2)+\lambda Y_1(\lambda R_1)]=0.
\end{aligned}
\label{e5.15}
\end{equation}
For $n=3$, $G$ has the form
\begin{equation}
\begin{aligned}
&G(r,\xi, t)\\
&=\frac{2\xi}{r}\displaystyle {\sum_{n=1} ^{\infty}} \frac{(b_2 ^2 +R_2 ^2 \lambda_n ^2)\Psi_n (r)\Psi_n (\xi)e^{-\frac{\epsilon}{2} \lambda_n ^2 t}}
{(R_2 -R_1)(b_1 ^2 +R_1 ^2 \lambda_n ^2)(b_2 ^2 +R_2 ^2 \lambda_n ^2)+(b_1 R_2 +b_2 R_1)(b_1 b_2 +R_1 R_2 \lambda_n ^2)},
\end{aligned}
\label{e5.16}
\end{equation}
where
\begin{equation}
\begin{aligned}
\Psi_n (r)&=b_1 \sin[\lambda_n (r-R_1)]+R_1 \lambda_n cos[\lambda_n (r-R_1)]\\
& b_1=-(q_1/\epsilon) R_1+1,\,\,\, b_2 =(q_2/\epsilon) R_2 -1,
\label{e5.17}
\end{aligned}
\end{equation}
and $\lambda_n$  are positive solutions of the transcendental  equation 
\begin{equation}
(b_1 b_2 -R_1 R_2 \lambda^2)\sin [\lambda(R_2 -R_1)]+\lambda(R_1 b_2 +R_2 b_1)\cos[\lambda(R_2 -R_1)]=0.
\label{e5.18}
\end{equation}
\end{theorem}

\begin{proof}
From theorem 2.1, we get $u=-\epsilon \frac{a_r}{a}$ where $a$ satisfies
\begin{equation}
\begin{aligned}
&a_t = \frac{\epsilon}{2}[a_{rr} +\frac{(n-1)}{r}a_r],\,\,\, R_1<r<R_2,\,\,\,t>0\\
&a(r,0)=e^{-\frac{\int_0^ r q(s) ds}{\epsilon}}, \,\,\,R_1<r<R_2\\
&\epsilon a_r(R_i ,t)+q_i a(R_i ,t)=0, \,\,\,i=1,2, \,\,\,t>0.\\
\end{aligned}
\label{e5.19} 
\end{equation}
Then $G$ given by \eqref{e5.13}-\eqref{e5.15} and \eqref{e5.16}-\eqref{e5.18} are just the Green's 
functions for the boundary value problem \eqref{e5.19} for $n=2$ and $n=3$ respectively, see \cite{p1}.  
Therefore the formula \eqref{e5.12} for $u$ follows. 
The formula for $\rho$ is obtained by the method of characteristics in the same way as in the previous
theorem and the details are omitted. 
\end{proof}

\subsection{Remarks on large time behaviour}
The large time behaviour of the velocity follows easily from the explicit formulae. In Case $1$, where we considered the problem in
the domain $\{(x,t) : |x|<R, t>0\}$, $u$ given by \eqref{e5.4}-\eqref{e5.6} has the following asymptotic form :

\begin{equation*}
\lim_{t \rightarrow \infty}u(x,t) = -\epsilon \frac{x}{|x|}\frac{\int_0^R \partial_r{G(r,\xi)}e^{-\frac{1}{\epsilon}\int_0 ^\xi q_0 (s)ds} d\xi}
{\int_0^R G(r,\xi) e^{-\frac{1}{\epsilon}\int_0 ^\xi q_0 (s)ds} d\xi}
\end{equation*} 
where for $n=2$,
\begin{equation*}
G(r,\xi)=\xi J_0 (\frac{\mu _1 \xi}{R})J_0 (\frac{\mu _1 r}{R}),
\end{equation*}
$\mu_1$ being the first positive solution of the equation 
$\mu J_1 (\mu)-kR J_0(\mu)=0$. For $n=3$,
\begin{equation*}
G(r,\xi)=\frac{\xi}{r} \sin (\frac{\mu_1 r} {R}) \sin (\frac{\mu_1 \xi} {R}),
\end{equation*}
$\mu_1$ being the first positive solution of the equation 
$\mu \cot(\mu)+\frac{q_B}{\epsilon}-1=0$.\\
Similarly, for Case $2$ where we consider the problem in the domain $\{(x,t) : R_1<|x|<R_2, t>0\}$, $u$ given by \eqref{e5.12}-\eqref{e5.18}
has the following asymptotic form :
\begin{equation*}
\begin{aligned}
 \lim_{t \rightarrow \infty}u(x,t)=-\epsilon \frac{x}{|x|}\frac{\int_{R_1}^{R_2} \partial_r{G(r,\xi)} e^{-\frac{1}{\epsilon}\int_0 ^\xi q_0 (s)ds}d\xi}
{\int_{R_1}^{R_2} G(r,\xi) e^{-\frac{1}{\epsilon}\int_0 ^\xi q_0 (s)ds} d\xi}
\end{aligned}
\end{equation*}
where for $n=2$,
\begin{equation*}
 G(r,\xi)= \xi H_1(r) H_1(\xi) 
\end{equation*}
where $H_1$ is given by \eqref{e5.14} with $\lambda_1$ being the first positive 
root of the equation \eqref{e5.15}. For $n=3$, 
\begin{equation*}
G(r,\xi)=\frac{\xi}{r} \Psi_1 (r)\Psi_1 (\xi)
\end{equation*}
where 
$\Psi_1 (r)=b_1 \sin[\lambda_1 (r-R_1)]+R_1 \lambda_1\cos[\lambda_n (r-R_1)]$, $b_1=-(q_1/\epsilon) R_1+1$, $b_2 =(q_2/\epsilon) R_2 -1$,
and $\lambda_1$ is the first positive solution of the equation \eqref{e5.17}.

The asymptotic behaviour of $\rho$ is more complex and depends on whether mass flows out or flows in.
Let us consider the mass $m(t) = \int_{D}\rho(x,t) dx$. When $\rho$ is spherically symmetric and $D=\{x :|x|<R\}$, 
\begin{equation*}
m(t)=\omega_{n-1} \int_{0}^{R}r^{(n-1)}\rho(r,t)dr.
\end{equation*}
From the second equation of \eqref{e2.9}, we have
\begin{equation*}
\frac{dm(t)}{dt}=-q(R,t)p(R,t)
\end{equation*}
In the case when $\rho$ is spherically symmetric and $D=R_1<|x|<R_2$,  
\begin{equation*}
\frac{dm(t)}{dt}=-(q(R_2,t)p(R_2,t)-q(R_1,t)p(R_1,t)).
\end{equation*}
These show that the mass is conserved if velocity is zero on the space boundary.

\section{Rankine-Hugoniot conditions for the one-dimensional case}
Let us consider the
initial value problem 
\begin{equation}
\begin{aligned}
u_t + f(u)_x &=0,\\
v_t + (f'(u) v)_x &=0
\end{aligned}
\label{e6.1}
\end{equation}
with initial conditions
\begin{equation}
    u(x,0)=u_0(x),\,\,\,v(x,0)=v_0(x),
\label{e6.2}
\end{equation}
with the flux satisfying  $f''(u)>0$, and $\lim_{|u|\rightarrow 0} \frac{f(u)}{u} =\infty$. The convex dual $f^{*}$ of $f$ is well defined
on $R^1$. There are different notions of solutions when the solution contains a measure component. Our aim in this section is to relate the notion
of solution introduced by LeFloch \cite{lef1} and the weak asymptotic solution of \cite{a1,s1}, for the one dimensional case.
 Let $u_0$ and $v_0$ are of bounded variation, using Lax's formula, it was shown in \cite{lef1} that the pair $(u,v)$ given by
\begin{equation}
u(x,t) = (f^{*})'(\frac{x-y(x,t)}{t}),\,\,\,
 v(x,t)=\partial_x(\int_0^{y(x,t)} v_0(z) dz) ,
\end{equation}
where $y(x,t)$ is a minimizer of $\min_{y}\{\int_0^y u_0(z) dz + tf^{*}(\frac{x-y}{t})\}$
is a weak solution of \eqref{e6.1}, satisfying the initial condition \eqref{e6.2}.

Here we show that this solution satisfies the Rankine-Hugoniot conditions derived in \cite{a1,s1}.
\begin{theorem}\label{thm6.1}  Assume that the solution is smooth away except along a
curve $x=s(t)$ in a nbd of a point $(s(t_0),t_0)$. Then $v$ has the form $v(x,t) =v_{l}(x,t) +H(x-s(t)) v_{r}(x,t) +e(t) \delta_{x=s(t)}$,
and the solution satisfies the Rankine-Hugoniot condition 
\begin{equation}
\begin{aligned}
s'(t)&=\frac{f(u(s(t)+,t))-f(u(s(t)-,t))}{u(s(t)+,t)-u(s(t)-,t)},\\
e'(t)&= -s'(t)[v](s(t),t) +[v f'(u)](s(t),t).
\end{aligned}
\end{equation}
\end{theorem}
\begin{proof}
 Assume that $x=s(t)$ is a discontinuity
curve for $u$ and let $(s(t_0),t_0)$ be a point on it. Also assume that for a small nbd of this point, this curve is the
only discontinuity; that is except along this curve, $u$ and hence $v$ is smooth. Let $\phi$ be any test function supported in
this neighborhood. Now from the formula for $v$, 
\[
 v(x,t) =v_{l}(x,t) +H(x-s(t)) v_{r}(x,t) +e(t) \delta_{x=s(t)}
\]
The equation is understood  
\[
 v_t + (\bar{f}'(u)v)_x =0
\]
in the sense of distribution, where 
\[
\bar{f}'(u)=\int_0^1 f'(u(s(t)-,t)+\tau(u(s(t)+,t)-u(s(t)-,t))) d\tau.
\]
This means that
\[
\begin{aligned}
 \int_{x<s(t)}(v_l(x,t) \phi_t &+(f'(u_l(x,t)) v_l(x,t)) \phi_x)dx dt\\
 +\int_{x>s(t)}(v_r(x,t) \phi_t &+(f'(u_r(x,t)) v_r(x,t)) \phi_x)dx dt\\
 &+ \int_0^\infty e(t) \phi_t(s(t),t)+\bar{f}'(u) e(t) \phi_x(s(t),t) dt =0
\end{aligned}
\]
Using
\[
\begin{aligned}
\int_0^1 f'(u(s(t)-,t)&+\tau(u(s(t)+,t)-u(s(t)-,t))) d\tau\\ &=\frac{\int_0^1 \frac{df}{d\tau}(u(s(t)-,t)+\tau(u(s(t)+,t)-u(s(t)-,t)))d\tau}{u(s(t)+,t)-u(s(t)-,t)}\\
& =\frac{f(u(s(t)+,t))-f(u(s(t)-,t))}{u(s(t)+,t)-u(s(t)-,t)}=s'(t)
\end{aligned}
\]
we get,
\[
\begin{aligned}
 \int_{x<s(t)}(v_l(x,t) \phi_t +(f'(u_l) v_l) \phi_x)dx dt &+\int_{x>s(t)}(v_r(x,t) \phi_t +(f'(u_r) v_r) \phi_x)dx dt\\
 &+\int_0^\infty e(t) \frac{d\phi(s(t),t)}{dt}=0.
\end{aligned}
\]
By integration by parts and arranging terms we get, for all test functions supported in a neighborhood of $(s(t_0),t_0), t_0>0$
\[
 \int_0^\infty (-s'(t)[v](s(t),t) +[v f'(u)](s(t),t)-e'(t)) \phi(s(t),t) dt =0.
\]
This gives the Rankine-Hugoniot condition for the delta shock.
\end{proof}
The case of our interest is when $f(q)=\frac{q^2}{2}$. The pair $(u,v)$ given by
 \begin{equation}
u(x,t) = \frac{x-y(x,t)}{t},\,\,\,
 v(x,t)=\partial_x(\int_0^{y(x,t)} v_0(z) dz) ,
\end{equation}
where $y(x,t)$ is a minimizer of $\min_{y}\{\int_0^y u_0(z) dz + \frac{(x-y)^2}{2t}\}$,
is a weak solution of \eqref{e6.1} with $f(q)=\frac{q^2}{2}$,
satisfying the initial condition \eqref{e6.2}. This formula is also derived in \cite{j6}, using vanishing viscosity method.
In this case, the previous theorem gives
\begin{theorem}
 Assume that the solution is smooth away except along a
curve $x=s(t)$ in a nbd of a point $(s(t_0),t_0)$. Then $v$ has the form $v(x,t) =v_{l}(x,t) +H(x-s(t)) v_{r}(x,t) +e(t) \delta_{x=s(t)}$,
and the solution satisfies the Rankine-Hugoniot condition 
\begin{equation}
s'(t)=\frac{[\frac{u^2}{2}]}{[u]},\,\, e'(t)= -s'(t)[v](s(t),t) +[uv](s(t),t).
\end{equation}
\end{theorem}


\section{Radially symmetric solutions with prescribed mass condition and normal velocity at the origin}
The aim of this section is to find explicit global radial solution of
\begin{equation}
u_t+(u.\nabla)u=0,\,\,\,\rho_t +\nabla.(\rho u)=0,
\label{e7.1}
\end{equation}
with initial conditions
\begin{equation}
 u(x,0)=\frac{x}{r} q_0(r),\,\,\,\rho(x,0)=\rho_0(r)=r^{-(n-1)} p_0(r), \,\,\,r=|x|,
\label{e7.2}
\end{equation}
a condition on the mass 
\begin{equation}
 \int_{R^n} \rho(x,t) dx = p_B(t)
\label{e7.3}
\end{equation}
and prescribed normal velocity at the origin,
\begin{equation}
 \lim_{x\rightarrow 0}\frac{x}{r}u(x,t)= q_B(t).
\label{e7.4}
\end{equation}
The functions $q_0(r)$ and $p_0(r)$ are considered to be in the space of bounded variations and the functions $p_B, q_B$ are considered continuous. In fact, using a change of variables as in the proof of Therorem $4.2$
it is sufficient that $p_B$ is a constant which again corresponds to the conservation of mass condition. We would also like to mention that the vanishing viscosity solution for $u$ in the radial case corresponds to the case when $q_B(t)=0$.
We start with the following remark which gives us an idea regarding the form of the solution that we can expect.\\

{\bf{Remark}}:
If we consider the radial solutions of of the inviscid system \eqref{e1.2} of the form $(u,\rho)=\frac{x}{r} q(r,t),\rho(r,t)$, $r=|x|$ ,
then from Theorem $2.2$ it follows that $q$ and  $\rho$ satisfy the equations 
\begin{equation}
\begin{aligned}
&q_t +q q_r =0,\\
&\rho_t +\frac{1}{r^{(n-1)}}(r^{(n-1)} \rho q)_r =0.
\end{aligned} 
\notag
\end{equation}
The transformation $p=r^{(n-1)} \rho$ changes the above system to 
\begin{equation}
\begin{aligned}
&q_t +q q_r =0,\\
&p_t + (pq)_r =0, \,\,\, r>0,\,\,\,t>0.
\end{aligned}
\label{e7.5} 
\end{equation}
The initial value problem for this equation is well understood. This motivates us to seek a solution with initial data 
\begin{equation}
 u(x,0)=\frac{x}{r} q_0(r),\ \rho(x,0)=\rho_0(r),\,\,\,r=|x|
\notag
\end{equation}
The corresponding initial data for $(q,p)$ is then given by
\begin{equation}
 q(r,0)=q_0(r),\,\,\,p(r,0)=p_0(r)=r^{n-1}\rho_0(r)
\label{e7.6}
\end{equation}
Entropy weak solution of \eqref{e7.5} and \eqref{e7.6} is
\begin{equation}
 p(r,t)=\frac{r-r_0(r,t)}{t},\,\,\, q(r,t)=\partial_r({\int_0^{r_0(r,t)}\tau^{n-1}\rho_0(\tau) d\tau})
\notag
\end{equation}
where
$r_0(r,t)$ is a minimizer of
\begin{equation}
\min_{r_0\geq 0}\{\int_0^{r_0}q_0(\tau) d\tau+\frac{(r-r_0)^2}{2t}\}.
\notag
\end{equation}

This suggests the following formula for the solution of \eqref{e7.1} with initial conditions
\begin{equation}
 u(x,0)=\frac{x}{r}q_0(r),\,\,\, \rho(x,0)=\rho_0(r)
\notag
\end{equation}
namely,
\begin{equation}
 u(x,t)=\frac{x}{r}\frac{(r-r_0(r,t))}{t},\,\,\, q(r,t)=r^{-(n-1)} \partial_r({\int_0^{r_0(r,t)}\tau^{n-1}\rho_0(\tau) d\tau}),\,\,\,r=|x|
\notag
\end{equation}
But since the singularity at $r=0$ has to be taken care of, we need to impose appropriate conditions at the origin. We would like to note here
that the formula for $u(x,t)$ above was derived as the vanishing viscosity limit of the velocity component (see Theorem $4.3$).\\

Getting back to our original problem under discussion, we note that for initial data $u_0(x)$, in the space of bounded variation, the vanishing viscosity limit 
of the velocity component $u$ remains in the space of bounded variation, see \cite{j3}. The density $\rho$ is generally a measure. 
A local study of the Cauchy problem for \eqref{e7.1} and propagation of delta wave front was carried out by Albeverio and Shelkovich \cite{a1} 
in the frame work of the weak asymptotic method.

For the radial case we give here a global construction of $(u,\rho)$.  We also do not require
entropy condition on the initial data for $u$, as the equation is reduced to one dimensional Burgers equation
and hence we use the well-known one dimensional theory.
 As a first step of construction of $(u,\rho)$, we derive the equation for the radial case. 
 From Theorem $2.2$ it follows that $q$ and $p=r^{(n-1)}\rho$ satisfy the equations
\begin{equation}
 q_t+q q_r =0,\,\,\,p_t +(pq)_r =0,
\label{e7.7}
\end{equation}
with initial conditions
\begin{equation}
 q(x,0)=q_0(r),\,\,\, p(x,0)=p_0(r), \,\,\,r=|x|,
\label{e7.8}
\end{equation}
and boundary condition for $q$ and integral condition on $p$:  
\begin{equation}
 q(0,t) = q_B(t),\,\,\, \omega_{n-1}\int_0^\infty p(r,t) dr=p_B(t)
\label{e7.9}
\end{equation}
which is to be understood in a weak sense \cite{b1}. Weak formulation for the boundary condition for
$q$ is
\begin{equation}
 \begin{gathered}
  \text{either}\,\,\, q(0+,t)=q_B(t)\\
  \text{or} \,\,\,q(0+,t)\,\,\, \leq 0 \,\,\,\text{and} \,\,\,q^2(0+,t)\leq \,\,\,{q_B^{+}(t)}^2
 \end{gathered}
\label{e7.10}
\end{equation}

and for $p$ is
\begin{equation}
 \begin{gathered}
  \text{if} \,\,\,q(0+,t)>0\,\,\,
  \text{then} \,\,\,\omega_{n-1}\int_0^\infty p(r,t)dr=p_B(t).
 \end{gathered}
\label{e7.11}
\end{equation}

To describe the solution, we follow \cite{j1,j4,j5}. We introduce a class of paths in the
quarter plane $D=\{ (z,s) : z\geq 0, s \geq 0\}$.
For each fixed $(r,r_0,t), r \geq 0, r_0 \geq 0, t > 0$,
$C(r,r_0,t)$ denotes the following class of paths $\beta$.
 Each path is connected from the 
point
$(r_0,0)$ to $(r,t)$ and is of the form $z=\beta(s)$, where $\beta$ is a 
piecewise linear function of maximum three lines. On $C(r,r_0,t)$, we 
define a functional
\begin{equation}
J(\beta) = -\frac{1}{2}\int_{\{s:\beta(s)=0\}} (q_B(s)^{+})^2 ds 
+ \frac{1}{2} \int_{\{s:\beta(s) \neq 0\}} (\frac{d\beta(s)}{ds}^2)ds.
\label{e7.12}
\end{equation}

We call $\beta_0$ the straight line path connecting $(r_0,0)$ and $(r,t)$ 
which does not touch the space boundary $x=0$, namely $\{(0,t), 
t>0\}$. Then let
\begin{equation}
 A(r,r_0,t)= J(\beta_0) = \frac{(r-r_0)^2}{t}.
\label{e7.13}
\end{equation} 
Any $\beta \in C^{*}(r,r_0,t) = C(r,r_0,t)-{\beta_0}$
is made up of three pieces, namely lines
connecting $(r_0,0)$ to $(0,t_1)$ in the interior and $(0,t_1)$ to 
$(0,t_2)$ on the boundary and $(0,t_2)$ to $(r,t)$ in the interior.
For such curves, it can be easily seen from \eqref{e7.12} that
\begin{equation}
J(\beta) = J(r,r_0,t,t_1,t_2) = 
-\int_{t_1}^{t_2}\frac{(q_B(s)^{+})^2}{2}ds + 
\frac{r_0^2}{2 t_1} + \frac{r^2}{2(t-t_2)}.
\label{e7.14}
\end{equation}
For curves $\beta \in C^{*}(r,r_0,t)$ made up of two straight lines with one piece lying on the boundary $r=0$, we can write down a similar
expression as in the earlier case. 
 
It was proved in \cite{j5} that there exists $\beta^{*} \in 
C^{*}(r,r_0,t)$ 
and corresponding $t_1(r,r_0,t)$, $t_2(r,r_0,t)$ so that
\begin{equation}
\begin{aligned}
B(r,r_0,t)&= \min \{J(\beta) :\beta \in C^{*}(r,r_0,t)\}\\
        &= \min \{J(r,r_0,t,t_1,t_2): \,\, 0\leq t_1 < t_2 < t\}\\
        &=J(r,r_0,t,t_1(r,r_0,t),t_2(r,r_0,t))
\end{aligned}
\label{e7.15}
\end{equation} 
is Lipschitz continuous. Further
\begin{equation}
\begin{aligned}
m(r,r_0,t)&= \min\{J(\beta) : \beta \in C(r,r_0,t)\}\\
        & = \min \{A(r,r_0,t),B(r,r_0,t)\}
\end{aligned}
\label{e7.16}
\end{equation}
and
\begin{equation}
Q(r,t)= \min \{m(r,r_0,t) + \int_0^r  q_0(s) ds, \,\, 0\leq r_0< \infty\}
\label{e7.17}
\end{equation}
are Lipschitz continuous function in their variables.
Further, minimum in \eqref{e7.17}
is attained at some value of $r_{0}\geq0$, which depends on $(r,t)$; we
call it $r_0(r,t)$. If $A(r,r_0(r,t),t)\leq B(r,r_0(r,t),t)$,
\begin{equation}
\begin{aligned}
Q(r,t)= \frac{(r-r_0(r,t))^2}{2t} +\int_0^{r_{0}(r,t)} q_0(s) ds
\end{aligned}
\label{e7.18}
\end{equation}
and if  $A(r,r_0(r,t),t)> B(r,r_0(r,t),t)$,
\begin{equation}
Q(x,t)=J(r,r_0(r,t),t,t_1(r,r_0(r,t),t),t_2(r,r_0(r,t),t))+\int_0^{r_0(r,t)} q_0(s) ds.
\label{e7.19}
\end{equation}

Here and henceforth $r_0(r,t)$ is a minimizer in \eqref{e7.18} and in 
the case of \eqref{e7.19}, $t_2(r,t)=t_2(r,r_0(r,t),t)$ and
$t_1(r,t)=t_1(r,r_0(r,t),t)$. 
With these notations, we have the following result.

\begin{theorem}\label{thm7.1}  With $r_0(r,t), A(r,r_0(r,t),t), B(r,r_0(r,t),t), t_1(r,t),t_2(r,t)$ as defined
above, define
\begin{equation}
u(x,t) = \frac{x}{r}\begin{cases}
 \frac{r-r_0(r,t)}{t},&\text{if } A(r,r_0(r,t),t)< B(r,r_0(r,t),t),\\
 \frac{r}{t-t_1(x,t)},&\text{if } A(r,r_0(r,t),t)> B(r,r_0(r,t),t),
\end{cases}
\label{e7.20}
\end{equation}
and
\begin{equation}
P(r,t) = \begin{cases}
 -\int_{r_0(r,t)}^{\infty} p_0(z)dz,&\text{if } A(r,r_0(r,t),t)< B(r,r_0(r,t),t),\\
  \omega_{n-1} p_B(t_2(x,t)),
  &\text{if } A(r,r_0(r,t),t)>B(r,r_0(r,t),t). 
\end{cases}
\label{e7.21}
\end{equation}
and set
\begin{equation}
\begin{gathered}
\rho(x,t)= \frac{\partial_{r}(P(r,t))}{r^{(n-1)}}.
\end{gathered}
\label{e7.22}
\end{equation}
Then the distribution $(u(x,t),\rho(x,t))$ given by 
\eqref{e7.20}-\eqref{e7.22}  
satisfies \eqref{e7.1} in the region of smoothness and the Rankine-Hugoniot conditions \eqref{e3.5}, along a discontinuity surface.
Further it satisfies the initial conditions 
\eqref{e7.2}, mass conditions \eqref{e7.3} and normal velocity at the origin \eqref{e7.4} in the weak sense
\eqref{e7.10}-\eqref{e7.11}.
\end{theorem}

\begin{proof}
 Explicit formula for the entropy weak solution $q$ satisfying \eqref{e7.8}-\eqref{e7.10} is given 
in \cite{j1,j4}. This formula involves only a finite dimensional minimization, namely  
\eqref{e7.15}-\eqref{e7.17} in three variables $(t_1,t_2,y)$. The formula takes the form
\begin{equation}
q(r,t) = \begin{cases}
 \frac{r-r_0(r,t)}{t},&\text{if } A(r,r_0(r,t),t)< B(r,r_0(r,t),t),\\
 \frac{r}{t-t_1(x,t)},&\text{if } A(r,r_0(r,t),t)> B(r,r_0(r,t),t).
\end{cases}
\label{e7.23}
\end{equation}
To get the component $p$ we consider the problem for 
\begin{equation}
 P(r,t)=-\int_{r}^\infty p(s,t) ds
\label{e7.24}
\end{equation}
so that
\begin{equation}
\begin{gathered}
p(r,t)= \partial_{r}(P(r,t)).
\end{gathered}
\label{e7.25}
\end{equation}
It is easy to see from \eqref{e7.7}, \eqref{e7.8} and \eqref{e7.11} that $P$
must satisfy
\begin{equation*}
 P_t +q P_x =0, 
\end{equation*}
with initial condition
\begin{equation*}
 P(r,0)=-\int_r^\infty p_0(s) ds
\end{equation*}
and boundary condition

\begin{equation*}
 \begin{gathered}
  \text{if} \,\,\,q(0+,t)>0, \,\,\, 
  \text{then} \,\,\, P(r,t)dr=-\frac{1}{\omega_{n-1}}p_B(t).
 \end{gathered}
\end{equation*}
An explicit formula for the solution of this problem is given in \cite{j4} and has the form \eqref{e7.21}.
The formula \eqref{e7.20}-\eqref{e7.22} is then obtained from the transformation \eqref{e2.1} and
\eqref{e7.23}-\eqref{e7.25}. 

Now we show that $(u,\rho)$ given by \eqref{e7.20}-\eqref{e7.22} is a solution of \eqref{e7.1} in the
region of smoothness. This fact follows easily  as $(q,p)$
satisfies the system \eqref{e7.7} and the weak form of initial and boundary conditions \eqref{e7.8}-\eqref{e7.11} (see \cite{j4}). Also it was  
shown in \cite{j4}, that $(q,p)$ is a weak solution of \eqref{e7.7} and that the surface of discontinuity of
$(u,\rho)$ and $(q,p)$ are same. In order to show that
$(u,\rho)$ satisfies the Rankine-Hugoniot conditions, we use the corresponding 
conditions for the system \eqref{e7.7} from Theorem $6.2$.

Consider a solution $(q(r,t),p(r,t))$ of the system \eqref{e7.7} with a discontinuity on the surface
 $S(r,t)=0$, with $S(r,t)=r-s(t)$.
In a neighbourhood of this surface, we assume $(q,p)$ has the form
\begin{equation}
 q(r,t)=\bar{q}(x,t),\,\,\, p(r,t)=\bar{p}(r,t)+e(t)\delta_{r=s(t)}.
\label{e7.26}
\end{equation}
where $\bar{q},\bar{p}$ are smooth except on the surface $r=s(t)$, and $e(t)$ differentiable function of $t$.
The Rankine-Hugoniot condition for $(q,p)$ takes the form 
\begin{equation}
 \frac{ds(t)}{dt}=\frac{q_{+}+q_{-}}{2},\,\,\, \frac{de}{dt}=[qp]-[p]\frac{ds}{dt}.
\label{e7.27}
\end{equation}
We show that the distributions \eqref{e7.20}-\eqref{e7.22} satisfy the Rankine-Hugoniot conditions \eqref{e5.6}. 
As $u=\frac{x}{r} q, r^{(n-1)}\rho = p$, the coefficient of the measure $\delta$ in the density $\rho$, 
$\hat{e}$ and that in $p$ are related by
\begin{equation} 
r^{(n-1)}\hat{e}=e.
\label{e7.28}
\end{equation}
 Since $S(x,t)=S(r,t)=r-s(t)$, with $r=|x|$, we have
 \begin{equation*}
  S_t =-\frac{ds}{dt},\,\,\, \nabla_{x}S = \frac{x}{r}.
 \end{equation*}
The first equation of \eqref{e5.6} becomes
\begin{equation}
 S_t + u_{\delta}.\nabla_{x}S|_{\Gamma_t}=-\frac{ds}{dt}+\frac{q_{+}+q_{-}}{2}.\frac{x}{r}.\frac{x}{r}=-\frac{ds}{dt}+\frac{q_{+}+q_{-}}{2}=0,
\label{e7.29}
\end{equation}
where we used the first equation of \eqref{e7.27}. To verify the second equation in \eqref{e3.5}, we compute each terms in radial components.

\begin{equation}
\begin{aligned}
\frac{\delta \hat{e}}{\delta t}&=\frac{\partial \hat{e}}{\partial t} + G \frac{\partial \hat{e}}{\partial r}
=\frac{\partial \hat{e}}{\partial t} - \frac{S_t}{S_r} \frac{\partial \hat{e}}{\partial r}
=\frac{d\hat{e}}{dt}.
\end{aligned}
\label{e7.30}
\end{equation}
Also we know from \cite{a1},
\begin{equation}
\nabla_{\Gamma_t}.(\hat{e}u_\delta)=-2K G \hat{e},
\label{e7.31}
\end{equation}
where $K$ is the mean curvature of the surface of discontinuity,
$K=-\frac{1}{2}\nabla .\nu,\,\,\, \nu=\frac{x}{r}$
and $G=-\frac{S_t}{S_r}$. An easy computation gives
\begin{equation}
 K=-\frac{1}{2}\sum_{i=1}^n\frac{\partial {\nu_i}}{\partial {x_i}}
=-\frac{1}{2}\frac{n}{r}+\frac{1}{2}\sum_{i=1}^n\frac{x_i^2}{r^3}=-\frac{n-1}{2r}.
\label{e7.32}
\end{equation}

Using \eqref{e7.32} in \eqref{e7.31}, we get
\begin{equation}
 \nabla_{\Gamma_t}.(\hat{e} u_\delta)=\frac{n-1}{r}(-\frac{S_t}{S_r})\hat{e} =\frac{n-1}{r} \frac{dr}{dt} \hat{e}.
\label{e7.33}
\end{equation}
Finally,
\begin{equation}
\begin{aligned}
([\bar{\rho} u]-[\bar{\rho}] u_{\delta}]).\nabla_{x} S|_{\Gamma_t}&=(\frac{1}{r^{(n-1)}}.\frac{x}{r}[pq]-
\frac{1}{r^{(n-1)}}[p]\frac{(q_{+}+q_{-})}{2} \frac{x}{r})S_r.\frac{x}{r}\\
&=\frac{S_r}{r^{(n-1)}}([pq]-[p]\frac{ds}{dt}).
\end{aligned}
\label{e7.34}
\end{equation}

Since $S_r=1$, we get from \eqref{e7.30},\eqref{e7.33} and \eqref{e7.34},
\begin{equation*}
\begin{aligned}
\frac{\delta \hat{e}}{\delta t} + \nabla_{\Gamma_t}(\hat{e} u_{\delta})-([\rho u]-[\rho]u_{\delta}).\nabla_{x}S|_{\Gamma_t}
 &=\frac{d\hat{e}}{dt} + \frac{n-1}{r}\frac{dr}{dt} \hat{e}- \frac{S_r}{r^{(n-1)}}([pq]-[p]\frac{dr}{dt})\\
&=\frac{1}{r^{(n-1)}}[\frac{d(r^{(n-1)}\hat{e})}{dt}-[pq]-[p]\frac{dr}{dt}]\\
&=\frac{1}{r^{(n-1)}}[\frac{d{e}}{dt}-[pq]-[p]\frac{dr}{dt}] =0,
\end{aligned}
\end{equation*}
where in the last equality we used the second equation of \eqref{e7.27}. This completes the proof of the theorem.
\end{proof}









\end{document}